\documentclass[10pt]{amsart}
\usepackage{yfonts,amsthm}           
\usepackage{graphicx, float}
\usepackage{amssymb,enumerate}
\usepackage{epstopdf}
\newtheorem{theorem}{Theorem}

\newtheorem{lemma}[theorem]{Lemma}
\newtheorem{conjecture}[theorem]{Conjecture}

\title[The approximate Loebl--Koml\'os--S\'os conjecture]{The approximate Loebl--Koml\'os--S\'os conjecture and embedding trees in sparse graphs}
 \author[J. Hladk\'y]{Jan Hladk\'y}
 \address{Institute of Mathematics, Czech Academy of Science. \v Zitn\'a 25, 110 00, Praha, Czech Republic. The Institute of Mathematics of the Czech Academy of Sciences is supported by RVO:67985840.}
 \email{honzahladky@gmail.com}
\author[D. Piguet]{Diana Piguet}
\address{Institute of Computer Science, Czech Academy of Sciences, Pod Vod\'arenskou v\v e\v z\'i 2, 182~07 Prague, Czech Republic. With institutional support RVO:67985807}
\email{piguet@cs.cas.cz}
\author[M. Simonovits]{Mikl\'os Simonovits}
\address{R\'enyi  Institute, Budapest, Hungary}
\email{miki@renyi.hu}
\author[M. Stein]{Maya Stein}
\address{Centro de Modelamiento Matem\'atico,
Universidad de Chile, Beauchef 851,  Santiago Centro, RM, Chile}
\email{mstein@dim.uchile.cl}
\author[E. Szemer\'edi]{Endre Szemer\'edi}
\address{Department of Mathematics, Rutgers University, 110 Frelinghuysen Rd., Piscataway, NJ~08854-8019, USA}


\begin{document}
\maketitle
\begin{abstract}
Loebl, Koml\'os and S\'os conjectured that every $n$-vertex graph $G$ with at least~$n/2$ vertices
of degree at least~$k$ contains each tree $T$ of order $k+1$ as a
subgraph. We give a sketch of a proof of the approximate version of this conjecture for large values of $k$. 

For our proof, we use a structural decomposition which can be seen as an analogue of Szemer\'edi's regularity lemma for possibly very sparse graphs. With this tool, each graph can be decomposed into four parts: a set of vertices of huge degree, regular pairs (in the sense of the regularity lemma), and two other objects each exhibiting certain expansion properties. We then exploit the properties of each of the parts of~$G$ to embed a given tree $T$.

The purpose of this note is to highlight the key steps of our proof. Details can be found in [arXiv:1211.3050].
\end{abstract}

%\pagebreak

\section{Introduction}
Szemer\'edi's Regularity Lemma from 1975 allows to decompose each dense graph into a bounded collection of random-like subgraphs. The lemma and its variants have found numerous applications in graph theory, number theory, and theoretical computer science. In particular, it is crucial to some developments on extremal problems concerning dense graphs, in the last two decades. On the other hand, extremal problems concerning sparse graphs have been lacking a general framework. We present a new tool which generalizes the Regularity Lemma and which applies to all graphs. This tool seems particularly suitable for embedding trees, even in very sparse graphs. As an application, we prove an approximate version of the Loebl--Koml\'os--S\'os Conjecture.

The Loebl--Koml\'os--S\'os Conjecture is a typical example of a problem in extremal graph theory. These often are of the following type: Does a certain density condition imposed on a graph of order $n$ guarantee a given subgraph? Statements of this spirit include Mantel's theorem, which states that an average degree of more than $n/2$ ensures a triangle as a subgraph, the more general Tur\'an theorem, and Dirac's theorem, which states that a minimum degree of at least $n/2$  forces the appearance of a Hamilton cycle. Other prominent results include the Erd\H os--Stone--Simonovits theorem which determines asymptotically the average degree threshold for appearance of any fixed graph, or the solution~\cite{PoSe} of the P\'osa--Seymour conjecture about containment of powers of Hamilton cycles. 

Some of the progress in the area in the last two decades has been enabled by the developments around the regularity lemma.  The regularity lemma allows to approximate an original graph by a so-called cluster graph. The point in doing so is that an original combinatorial problem translates to an easier one on the cluster graph. Let us illustrate this fundamental feature with the examples mentioned above: the modern, regularity lemma based approach reduces the Erd\H os--Stone--Simonovits theorem to Tur\'an's theorem on the cluster graph. In a similar spirit, the proof of the P\'osa--Seymour conjecture about the appearance of the $k$-th power of a Hamilton cycle can be reduced to an easier question of tiling with copies of $K_{k+1}$ on the cluster graph level, an answer to which is given by the Hajnal--Szemer\'edi theorem. These and other applications of the regularity lemma in extremal graph theory are surveyed in~\cite{regu2,KuOs}.

Containment of trees is a particularly important case to study, as trees constitute a relatively simple graph class.
An easy greedy embedding argument shows that each graph with a minimum degree of at least $k$ contains each tree with $k$ edges. A graph formed by a union of cliques of order~$k$ shows that this result is optimal.

Two important conjectures have been made as to how the minimum degree condition can be relaxed. The first of these is the famous Erd\H os--S\'os conjecture from 1963: 
\begin{conjecture}\label{conjES}
Every graph of average degree greater than $k-1$ contains all trees with $k$ edges as subgraphs.
\end{conjecture}
Conjecture~\ref{conjES} trivially holds for the containment of a star with $k$ leaves, and it is a classical result of Erd\H os and Gallai~\cite{ErGa} that it also holds for paths (for more history, see~\cite{FurSim}). Further partial results include~\cite{bradob,Haxell:TreeEmbeddings, sacwoz, woz}. A proof of the conjecture for large graphs has been announced by Ajtai, Koml\'os, Simonovits and Szemer\'edi~\cite{AKSS07+}.

Loebl, Koml\'os, and S\'os (see~\cite{EFLS95}) conjectured the same assertion holds when replacing the average degree condition with the median degree.
\begin{conjecture}\label{conjLKS}
Every graph of median degree at least $k$ contains all trees with $k$ edges as subgraphs.
\end{conjecture}
Previous work on Conjecture~\ref{conjLKS} includes solutions which use additional restrictions on the host graph~\cite{Sof00, Dob02}, or on the trees~\cite{BLW00, PS2}. Most notably,  Conjecture~\ref{conjLKS} has been solved for large dense graphs, i.e.~for $k$ linear in $n$, in~\cite{HladkyPiguet, Cooley08}, building on an approximate version given in~\cite{PS07+}.  For the exact value $k=n/2$, this had been achieved earlier in~\cite{AKS95, Z07+}.

\begin{figure}[t]
\centering 
\includegraphics[scale=0.7]{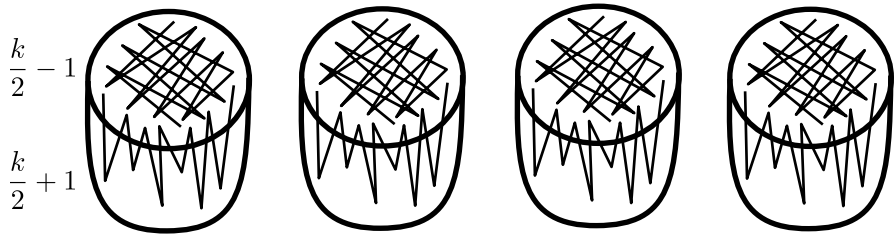}
\caption[Extremal graph for the Loebl--Koml\'os--S\'os Conjecture]{An extremal graph for the Loebl--Koml\'os--S\'os Conjecture.}
\label{fig:ExtremalGraph}
\end{figure}

Note that because of the stars, the median degree in Conjecture~\ref{conjLKS} has to be at least $k$. Further note that Conjecture~\ref{conjLKS} is almost best possible in the sense that we cannot decrease much the number of vertices, namely $n/2$, that are required to have degree at least $k$. 
For this, first assume that $n$ is even, and that $n=k+1$. Let $G^*$ be obtained from the complete graph on~$n$ vertices by deleting all edges inside a set of  $\frac
n2+1$ vertices. Then $G^*$ has $\frac n2-1$ vertices of degree $k$. It is easy to
check that $G^*$ does not contain the
path with $k$ edges (or  any other
tree with $k$ edges and independence number less than $\frac n2 +1$). Now, taking the union of 
several disjoint copies  of $G^*$ we obtain examples for other values of $n$. (And adding a small complete component we can get to {any} value of $n$.) See
Figure~\ref{fig:ExtremalGraph} for an illustration.

It is not difficult to see that each of the two conjectures implies that the Ramsey number of two trees $T_k$ and $T_m$ with $k$ and $m$ edges, respectively, is at most $k+m$. This is best possible for stars of even order, but not for all trees~\cite{GerencserGyarfas, HLT02}. 

Our main result is an approximate version of Conjecture~\ref{conjLKS}, which reads as follows.
\begin{theorem}\label{thmLKSappr}
For every $\varepsilon > 0$ there exists $k_0$ such that for every $k > k_0$, every $n$-vertex graph $G$ with at least $(1 +\varepsilon)n/2$ vertices of degree at least $(1+\varepsilon)k$ contains
each  tree $T$ with $k$ edges.
\end{theorem}

Previous results~\cite{AKS95,PS07+,Z07+,HladkyPiguet,Cooley08} on the dense case of Conjecture~\ref{conjLKS} relied on
Szemer\'edi's regularity lemma. The strategy of these proofs is explained in the next section. The (original) regularity lemma is void when the host graph is sparse, i.e., when $k=o(n)$. To circumvent this shortcoming, we present  an extension of the regularity lemma which is tailored to tree-embedding problems, and which applies even to sparse graphs. We then show how this decomposition, which we call \emph{sparse decomposition}, can be used to embed the tree $T$ given by Theorem~\ref{thmLKSappr}.

In this paper we show some of the key ideas behind the proof. The actual implementation of these ideas is technical and can be found in~\cite{LKSsparseapprox}, split into four parts~\cite{LKSsparseapprox1,LKSsparseapprox2,LKSsparseapprox3,LKSsparseapprox4} for publication purposes.

\section{The dense case}\label{sec:densecase}
In this section, we recall the solution of the dense approximate version of Conjecture~\ref{conjLKS} due to Piguet and Stein. Their proof provides several key ingredients which are common to the proof of Theorem~\ref{thmLKSappr}.
\begin{theorem}[\cite{PS07+}]\label{thm:PS}
For every $C,\varepsilon > 0$ there exists $k_0$ such that for every $k > k_0$ we have that every graph $G$ of order $n<Ck$ with at least $n/2$ vertices of degree at least $(1+\varepsilon )k$ contains
each tree $T$ with $k$ edges.
\end{theorem}
The proof follows a strategy typical for graph theory results which employ the regularity lemma. 
%(See~\cite{KuOs} for an extensive survey of applications of the regularity lemma in extremal graph theory.)
It has three main steps: partitioning~$T$, finding a suitable matching structure in the cluster graph~$\mathbf G$ of the graph~$G$, and embedding~$T$ into~$G$ using this matching structure.
The right way to partition~$T$ is given by the following lemma. We say that a subtree $F\subset T$ is \emph{adjacent} to a vertex $v\in V(T)\setminus V(F)$ if there is an edge from $v$ to a vertex in $F$.
\begin{lemma}\label{lemma_tree}
For each $\tau >0$, $k\in\mathbb N$, for any tree $T$ with $k$ edges there is a set $W=W_A\cup W_B\subset V(T)$, and a set $\mathcal T=\mathcal T_A\cup \mathcal T_B$ of disjoint subtrees covering all  of $V(T)-W$ such that
%For each $\tau >0$, $k\in\mathbb N$, for any tree $T$ with $k$ edges there is a set $\mathcal T=\mathcal T_A\cup \mathcal T_B$ of disjoint subtrees, and a set $W=W_A\cup W_B$ of vertices  such that $V(T)$ is the disjoint union of $W$ and $V(\bigcup \mathcal T)$, and
\begin{enumerate}[(a)]
\item the trees in $\mathcal T$ have order less than $\tau k$,
\item the trees in $\mathcal T$ are not adjacent to each other,
\item\label{e:W} $|W|<100/\tau$,
\item each side of the bipartition of $V(T)$ contains one of the sets $W_A$, $W_B$,
\item\label{E:TB} each tree in $\mathcal T_B$ is adjacent to only one vertex of $W$, and that vertex lies in $W_B$,
\item\label{E:TA} each tree in $\mathcal T_A$ is adjacent to at most two vertices of $W$, and these lie in $W_A$,
\item $|V(\bigcup \mathcal T_B)|< k/2$.\label{khalf}
\end{enumerate}
\end{lemma}
We call a tree in $\mathcal T$ \emph{internal} if it is adjacent to two vertices of $W$, and call it an \emph{end tree} otherwise.  Note that $\mathcal T_B$ only contains end trees by property~\eqref{E:TB}.

In order to obtain the partition from Lemma~\ref{lemma_tree}, we traverse $T$ from the leaves to a fixed root, sequentially chopping off `branches' of $T$ that have reached the critical size of $\tau k$. Similar strategies of dividing larger objects into smaller pieces have been used in other proofs employing  embedding with the regularity method. As we shall see, the partition of $T$ from Lemma~\ref{lemma_tree} will be useful for the sparse case as well.

The bulk of the work is on finding a suitable structure in the graph $G$. To this end we use the regularity lemma~\cite{Sze78}. Let us first introduce the key notion of regular pairs.
Given $\eta > 0$, a pair $(A,B)$ of disjoint sets is \emph{$\eta$-regular} if $|d(A,B) - d(U,W)| < \eta$ for each $U \subseteq A, W \subseteq B$ with $|U| > \eta |A|, |W| > \eta |B|$. The regularity lemma then reads as follows.
\begin{lemma}[Regularity lemma]\label{lem:regu}
For each $\eta>0, m\in\mathbb N$ there are $n_0, M$ such that every graph on $n>n_0$ vertices allows for a partition of all but at most $\eta n$ of its vertices into $m<k<M$ sets (the `clusters') such that all but at most $\eta k^2$ pairs of clusters form $\eta$-regular pairs.
\end{lemma}

With the help of Lemma~\ref{lem:regu}, we regularize the graph $G$ and obtain a cluster graph $\mathbf G$ with clusters of size $\nu k$.
Let $\mathbf L$ be the set of those clusters of $\mathbf G$ whose typical vertices have degree more than~$k$. Piguet and Stein show that $\mathbf G$ contains a matching $M$ plus two adjacent vertices $A,B\in\mathbf L$ such that the vertices in~$A$ have degree more than $k$ into $M$, and those in $B$ have degree larger than $k/2$ into $V(M)\cup\mathbf L$.

The tree $T$ can be embedded into $G$ by suitably mapping the vertices of $W_A$ into $A$, the vertices of $W_B$ into $B$, and packing subtrees $\mathcal T$ (viewed as bipartite graphs) into the edges of $M$ and the edges $E_\mathbf{L}$ emanating from $\mathbf L$ using basic properties of regular pairs. Here, it is crucial we choose the parameter $\tau$ for  Lemma~\ref{lemma_tree} such that $\tau\ll \nu$. That is, individual subtrees of $\mathcal T$ are much smaller than the clusters.

The large degrees of $A$ and $B$ into $V(M)\cup\mathbf L$ guarantee that there is enough space for embedding all trees from $\mathcal T$. 
More precisely, each time we wish to embed a tree $T'\in \mathcal{T}$, there are two things we have to ensure. The first is that  we have enough free space in the neighbourhood of either $A$ or $B$ to map the root of $T'$. The second is that we have sufficient free space in some regular pair meeting this neighbourhood, to map the rest of $T'$. For this, a degree of about $k$ for the vertices in $A$ and a degree of about $k/2$ for the vertices in $B$ is sufficient.

\section{The sparse decomposition}\label{sec:sparsedecomposition}
In this section we introduce the basis of our  proof of Theorem~\ref{thmLKSappr}, the \emph{sparse decomposition}. This tool has been conceived by Ajtai, Koml\'os, Simonovits and Szemer\'edi during their work on the Erd\H os--S\'os conjecture. It allows to decompose any given graph, after a removal of a small fraction of the edges, into four sets:  a set $\Psi$ of vertices of high degree, a graph  $\mathbf{G}_\mathrm{reg}$ whose edges that span regular pairs, an expanding graph  $G_\mathrm{exp}$, and a set $\mathfrak{A}$ of vertices which has a yet different expansion property. 

Throughout this section, let us fix a graph $G$ on $n$ vertices, and let~$k$ have the same order of magnitude as the average degree of $G$.\footnote{This setting is compatible with the one of Theorem~\ref{thmLKSappr}. Indeed, a straightforward calculation gives that in that case the average degree $d$ of $G$ satisfies $d>k/2$. If, on the other hand, $d\ge 2k$ then there is no need to use the sparse decomposition as we can pass to a subgraph with minimum degree at least $d/2\ge k$, and embed $T$ greedily.}
We use Greek majuscules, and minuscules to denote sufficiently large, and sufficiently small constants, respectively. Some of these constants may depend on $G$ and $k$, but are absolutely bounded from above and from below. We make the subtle relations between these constants explicit only when it adds to the clarity of this rough sketch.

The first step of the sparse decomposition is to separate the vertices of very high degree from those of comparatively low degree. By deleting only a few well-chosen edges, we arrive at a subgraph~$G'$ of $G$ that has a gap in its degree sequence. More precisely, there are numbers $\Omega^*$ and $\Omega^{**}$ with $\Omega^*\ll\Omega^{**}$ such that no vertex of $G'$ has degree between $\Omega^*k$ and $\Omega^{**}k$ in $G'$. Let us indicate how to create the gap. We fix constants $1\ll\Omega_1\ll \Omega_2\ll\ldots \ll\Omega_{\lceil \frac 1 \varepsilon\rceil+1}:=\infty$. There is an index $i\in[\lceil \frac 1 \varepsilon\rceil+1]$ such the total degree of the vertices $v$ with $\deg(v)\in [\Omega_i k,\Omega_{i+1}k)$ is at most $\varepsilon kn$. Deleting the edges incident with these vertices almost yields the gap with $\Omega^*=\Omega_i$, and $\Omega^{**}=\Omega_{i+1}$. The problem is that the edge deletion may cause degrees of other vertices fall into the forbidden region $[\Omega_i k,\Omega_{i+1}k)$. This can be resolved using an additional argument, which we omit here. 

Let $\Psi$ denote the set of all vertices of degree at least $\Omega^{**} k$. The mere structural information about the vertices $\Psi$ is that they have huge degrees. On the other hand, this property turns out to be so powerful for tree embeddings  that it compensates the lack of any finer description.

Before proceeding with the decomposition, we need a few concepts. The \emph{density} of a  bipartite graph $D = (U,W; F)$  is $d(D):=\frac{| F| }{| U| | W| }$, where $F=E(D)$ are the edges of $D$. An \emph{$(m, \gamma)$-dense spot} in a graph
 is a non-empty bipartite subgraph $D$  with density $d(D)> \gamma$ and minimum degree $\delta (D) > m$. A graph $H$ is \emph{$(m, \gamma)$-nowhere-dense} if it does not contain any $(m, \gamma)$-dense spot. 
 
Let $\mathcal D$ be a maximal set of edge-disjoint $(\gamma k, \gamma)$-dense spots in $G'-\Psi$. Let $G_\mathrm{exp}$ be the $(\gamma k, \gamma)$-nowhere dense graph obtained from $G'-\Psi$ by removing the edges of $\mathcal D$. (We chose the name $G_\mathrm{exp}$ for this graph in order to emphasize its expansion property given by the fact it is nowhere dense.) We now  sequentially remove from $G_\mathrm{exp}$ any vertex of degree less than $\rho k$, where $\rho$ is a certain constant much smaller than $\varepsilon$, but greater than $\gamma\ll \rho\ll \varepsilon$. Note that in the cleaning procedure we lose less than $\rho kn$ edges, and the obtained graph, which we still call $G_\mathrm{exp}$, has minimum degree at least~$\rho k$.

The next step consists of regularizing the dense spots in $\mathcal D$. By this we mean we wish to find a graph spanning almost all of $\bigcup\mathcal D$, and consisting of clusters that pairwise mostly form highly regular pairs, in the sense above. For each single one of these spots this is possible by Lemma~\ref{lem:regu} above, but 
such a naive regularization of each dense spot separately is useless. Indeed, for embedding $T$ we may need to traverse many different spots $\mathcal D$. Thus the cluster structure of different dense spots must agree on their intersection. 

So, consider all the Venn cells $\mathcal V$ with respect to the system $\{U,W: (U,W; F)\in\mathcal D\}$. We shall not attempt to regularize those Venn cells $\mathcal A\subset \mathcal V$ which are of size less than $\alpha k$ (for $\alpha\ll\eta$), as those cells themselves may be smaller than the anticipated cluster sizes. We now construct an auxiliary graph $\mathcal G=(\mathcal V\setminus \mathcal A,\mathcal E)$ on the larger Venn cells,  joining two Venn cells $X,Y$ with an edge if there is a dense spot $(U,W;F)\in\mathcal D$ with $X\subset U, Y\subset W$. 

Let us sketch how to regularize simultaneously all the dense pairs corresponding to the edges of $\mathcal G$ in this setting. It can be shown that the maximum degree of $\mathcal G$ is bounded from above by a constant $\Delta$ that is  independent of $k$. By Vizing's theorem, we can cover $\mathcal G$ with $\Delta+1$ matchings $M_1,\ldots,M_{\Delta+1}$. We follow the idea of Szemer\'edi's proof of the regularity lemma, pumping up the mean square energy when refining an irregular partition. The key difference is that we track $\Delta+1$ mean-square energies, one for each matching $M_i$, rather than just a single one. This is similar to the proof of the multi-coloured version of the regularity lemma, which tracks a mean-square energy for each colour separately. We thus obtain a system of clusters,  of size $\nu k$, say, refining $\mathcal V\setminus \mathcal A$ and regular pairs between some of these clusters. Let $\mathbf{G}_\mathrm{reg}$ be the graph spanned by the regular pairs of positive density.

It remains to make use of vertices in $\mathfrak{A}:=\bigcup \mathcal A$, i.e., those in small Venn cells. An elementary double-counting argument gives the following expansion property of $\mathfrak{A}$, which we call the \emph{$(\Lambda,\beta,\gamma)$-avoiding property}:  For every $X \subseteq V(G)$
with $|X| \leq \Lambda k$ for all but at most $\beta k$ vertices $v \in \mathfrak A$ there is a dense spot
$D\in \mathcal  D$ which contains $v$ and which satisfies $|X \cap V(D)| \leq\gamma^2 k$. The constant $\beta$ will be chosen much smaller than $\gamma$, but still larger than $\tau\ll\beta\ll\gamma$.

Putting all of the above together, we obtained a sparse decomposition $(\Psi,\mathbf G_\mathrm{reg},G_\mathrm{exp},\mathfrak{A})$ which captures all but at most $o(kn)$ edges of $G$.

\section{Embedding the tree $T$}\label{sec:embedding}
The proof of Theorem~\ref{thm:PS} as outlined in Section~\ref{sec:densecase} is a combination of two elements: a global embedding strategy based on the matching structure given by the clusters $A,B$, and the matching~$M$, and a local strategy applied sequentially for embedding the individual subtrees from $\mathcal T$. The local strategy there is  the standard technique of filling up regular pairs. 

Also in the proof of Theorem~\ref{thmLKSappr} we shall find a suitable global structure, now in the sparse decomposition instead of in the cluster graph. 
We will discuss this structure later on. Before, we indicate local strategies for embedding subtrees $\mathcal T$ in each of the ingredients  $\Psi$, $\mathbf{G}_\mathrm{reg}$, $G_\mathrm{exp}$, $\mathfrak{A}$ of  the sparse decomposition. The starting point of using either $\Psi$, $\mathbf{G}_\mathrm{reg}$, $G_\mathrm{exp}$, or $\mathfrak{A}$ is that in our sequential embedding procedure we wish to extend the partial embedding from a vertex with a substantial degree into the respective object. For this, suppose that $U\subset V(G)$ is the set of vertices used in earlier steps of the embedding. 

To work with  $\mathfrak{A}$, we use the avoiding property. Say we have to embed a tree $T'\in\mathcal T$. For simplicity, let us assume that $T'\in\mathcal T_B$; this guarantees that $T'$ is an end tree. Suppose its parent in $W_B$ has been embedded already in a vertex $v$ of degree more than $\beta k$ into~$\mathfrak{A}\setminus U$. Then, by the definition of the avoiding property, there is a neighbour of $v$ in~$\mathfrak{A}\setminus U$ that is contained in a dense spot $D$ which does not meet~$U$ much.
We can place the first vertex of $T'$ appropriately into~$D$. We then use the minimum degree $\gamma k$ of $D$ to embed the rest of $T'$ greedily. Here, we use that $\tau+\gamma^2\le\gamma$.

Next, we show how to use $G_\mathrm{exp}$. Again, suppose we are in the process of embedding a subtree $T'\in\mathcal T_B$. Say $x$ is the last vertex of $T'$ that has been embedded already, and we now wish to embed the children of $x$.  
Assume the image $v$ of $x$ has neighbourhood $N_v\subset V(G_\mathrm{exp})$ of size at least $\rho k/2$ in $G_\mathrm{exp}-U$. (Below it will become clear why we may assume this.) Since $G_\mathrm{exp}$ does not contain any $(\gamma k, \gamma)$-dense spots, in particular not between $U $ and $N_v$, we know that most vertices in $N_v$ have less than $\rho k/3$ neighbours in $U$. (Here, we used that $\gamma\ll\rho$.) Thus, it is possible to embed the children of $x$ in vertices that have degree at least $\rho k/2$ in $G_\mathrm{exp}-(U\cup x)$, placing them in equally good positions as~$x$ earlier. Following this strategy successively, we manage to embed all of $T'$.

Regular pairs in~$\mathbf{G}_\mathrm{reg}$ are used in the usual way for embedding trees of $\mathcal T$. That is, we view these trees as bipartite graphs, and embed them in regular pairs using the regularity property.

It only remains to explain the role of $\Psi$. 
%Observe that it is enough to prove Theorem~\ref{thmLKSappr} only for graph for which the set of vertices of degrees at least $(1+\epsilon)k$ is independent. In particular, we can assume that $\Psi$ is independent. 
This set is used very rarely for embedding; in particular we always have $|U\cap \Psi|<\lambda k/2$.
The idea is that after mapping a vertex $x\in T$ to a vertex $v\in\Psi$, we have an affluence of choices to extend the embedding. However, the huge degree of $v$ alone is not enough. For example, we cannot map non-leaf vertices of $T$ to leaves of $G$, and these may potentially comprise the entire neighbourhood of $v$. To circumvent this issue, we employ a rather delicate cleaning procedure prior to starting the embedding. That is, we find a set $\Psi'\subset\Psi$ of vertices which have degree at least $\Omega'k$ (for suitable $\Omega^*\ll\Omega'\ll\Omega^{**}$) into a `useful part of $G$' in such a way that we do not lose many edges during the cleaning. 
%DDD
Having done so, we wish to map the neighbours of~$x$ to neighbours of~$v$ that send not more than a few edges to~$U$. This will guarantee that we can extend the embedding avoiding~$U$ in subsequent steps.  To this end, consider the set $\tilde U=\{u\in V(G):\deg(u,U)\ge \lambda k\}$. We have $\tilde U\subset \{u\in V:\deg(u,U\setminus \Psi)> \lambda k/2\}$, and double-counting the edges between $\tilde U$ and $U\setminus \Psi$ gives
$$|\tilde U|< \frac{\Omega^*k|U|}{\lambda k/2}\ll \Omega' k\;.$$
In particular, the neighbours of $x$ can be embedded outside of $\tilde U$.

\medskip
In Lemma~\ref{lemma3}, we describe the structural counterpart of the matching structure $(A,B,M)$ from the dense case (again, the structure is much simplified for presentation reasons). Note that this global structure must combine properties of all the objects $\Psi$, $\mathbf G_\mathrm{reg}$, $G_\mathrm{exp}$, $\mathfrak{A}$ as it could happen that none of them alone suffices for embedding~$T$. 

We write $\mathbb L$ for the set of those vertices of $G$ that have degree at least $(1+\varepsilon)k$ in $G$. We call a collection $\mathcal M$ of $\eta$-regular pairs of positive density with clusters of sizes $\mu k$ an \emph{$(\eta,\mu)$-regular matching} if all the regular pairs are disjoint.
\begin{lemma}\label{lemma3}
The graph $G$ contains two disjoint sets $\mathbb{X},\mathbb{Y}\subseteq \mathbb L$, and an $(\eta,\mu)$-regular matching $\mathcal M$ with the following properties. 
  \begin{enumerate}[(i)]
 \item\label{e:mindeg} The bipartite graph $G[\mathbb{X},\mathbb{Y}]$ has minimum degree at least $100/\tau$.
 \item\label{e:X} The vertices in $\mathbb{X}$ each have degree at least $(1+\varepsilon/2)k$ into $$Q:=\big(\Psi\cup V(G_\mathrm{exp})\cup\mathfrak{A}\cup\mathbb L\cup  V(\mathcal M)\big)\setminus (\mathbb{X}\cup\mathbb{Y})\;,$$
  \item\label{e:Y} The vertices in $\mathbb{Y}$ have degrees at least $(1+\varepsilon/2)k/2$ into $Q$, and
  \item $\mathcal M$ and $\mathbb{X}\cup\mathbb{Y}$ are disjoint.
   \end{enumerate}
  \end{lemma}

The sets $\mathbb{X}$ and $\mathbb{Y}$ will host $W_A$ and $W_B$, that is, we can think of $\mathbb{X}$ and $\mathbb{Y}$ being counterparts to the sets $A$ and $B$ from the dense case. Property~\eqref{e:mindeg} then guarantees that the edges between $W_A$ and $W_B$ can be embedded greedily (cf.~Lemma~\ref{lemma_tree}\eqref{e:W}). Properties~\eqref{e:X} and~\eqref{e:Y} guarantee that subtrees~$\mathcal T_A$ and $\mathcal T_B$ can be embedded. Instead of embedding $\mathcal T_A$ and $\mathcal T_B$ using the matching $M$ and edges $E_\mathbf{L}$ as in the proof of Theorem~\ref{thm:PS}, we have to make use of embedding techniques developed above.

\medskip
There are two major issues with the indicated approach. The first difficulty is encountered when embedding an internal tree $T'\in \mathcal T_A$. As such a tree may be adjacent to two vertices of $W_A$, and we plan to embed $W_A$ in $\mathbb{X}$, we have to return to $A$ after embedding~$T'$.

 To understand this difficulty better, it is instructive to first see how an internal tree $T'$ with head $x\in W_A$ and tail $y\in W_A$ is embedded in the dense case (head and tail are the two vertices from Lemma~\ref{lemma_tree}\eqref{E:TA}). Say $x$ has been embedded into vertex $v\in A$. We choose an edge $XY\in M\cup E_\mathbf{L}$ that will host $T'$, such that  $X$ is an edge in the cluster graph. The regularity of $(A,X)$ guarantees that the embedding of $T'$ can be extended from $x$, but also, that after embedding~$T'$ we can embed~$y$ back in $A$. 
 
 In the sparse case, we do not have any similar property for the set $Q$. To resolve this issue, we introduce certain cleaning procedures which guarantee that the last vertex before a tail of an internal tree is always embedded in a vertex of $Q$ which has degree at least $100/\tau$ into $\mathbb{X}$.

\begin{figure}[htbp] \centering \includegraphics[scale=0.8]{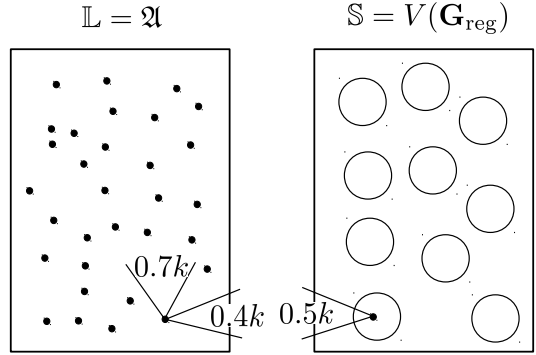}
\caption{The graph $G$ has a set $\mathbb{L}$ of  $\frac{5n}9$ vertices  of degree $1.1k$, and a set $\mathbb{S}$ of $\frac{4n}9$ vertices of degree $0.5k$. For each $v\in \mathbb{L}$, we have $\deg(v,\mathbb L)=0.7k$ and $\deg(v,\mathbb S)=0.4k$. For each $v\in \mathbb{S}$, we have $\deg(v,\mathbb L)=0.5k$ and $\deg(v,\mathbb S)=0$. Further, there is a sparse decomposition of $G$ such that the union of the dense spots $\mathcal D$ covers all the edges of $G$.  The dense spots $\mathcal D$ intersect in such a way that the small Venn cells $\mathcal A$ cover $\mathbb{L}$. 
The cluster graph $\mathbf{G}_\mathrm{reg}$ is edgeless, and all vertices have degree less than $k$ into $\Psi\cup V(G_\mathrm{exp})\cup\mathfrak{A}\cup\mathbb L$. Thus, for our set $\mathbb{X}$, we need to find $\mathcal M$ elsewhere. In this particular case, one can show that there is a regular matching between $\mathbb L$ and $\mathbb S$ covering almost all of $\mathbb S$, and that such a matching is a good choice for $\mathcal M$.}
\label{fig:whyEnhancing}
\end{figure}

The second difficulty arises when constructing the regular matching $\mathcal M$. In  analogy to the dense case, it would be a natural guess that $\mathcal M$ is a matching in $\mathbf G_\mathrm{reg}$. In Figure~1 we give an example that it is not always possible to choose $\mathcal M$ like this.  

Given the example of Figure~1, one might wonder why we bother to construct the cluster graph~$\mathbf{G}_\mathrm{reg}$ at all. The answer is that for constructing the regular matching $\mathcal M$, the graph $\mathbf{G}_\mathrm{reg}$ is of help, either directly, or via the information it gives by lacking a suitable matching. 

\section{Concluding remarks}
Let us conclude with several comments.
\begin{itemize}
\item Our proof builds on techniques developed by Ajtai, Koml\'os, Simonovits and Szemer\'edi for their work on the Erd\H os--S\'os conjecture. However, there is a substantial difference between the proofs already on the level of the sparse decomposition. In their proof, a suitable matching structure can always be found in the cluster graph $\mathbf{G}_\mathrm{reg}$. That means examples like the one in Figure~\ref{fig:whyEnhancing} do not enter the picture in the Erd\H os--S\'os conjecture. \\
\item Similarly as in the Erd\H os-S\'os setting, it seems that our approach can be combined with the stability approach of Simonovits. We hope to resolve the Loebl--Koml\'os--S\'os conjecture exactly, for $k$ sufficiently large (this is work in progress).\\
\item The sparse decomposition of a graph is not uniquely determined, and can actually vary vastly. This is caused by the arbitrariness in the choice of the dense spots from which we obtain the regularized graph $\mathbf{G}_\mathrm{reg}$.
This situation is in acute contrast with the situation of decomposition of
dense graphs (given by the regularity lemma).  Indeed, in
the dense setting the structure of the cluster graph is essentially
unique, cf.~\cite{AlShSt:DistinctyRegularityPartitions}.\footnote{In order to have uniqueness, the setting needs to be somewhat strengthened; see Theorem~1 and Theorem~2
in~\cite{AlShSt:DistinctyRegularityPartitions}. The uniqueness phenomenon can be nicely expressed in the language of graph limits~\cite{borgs-2008}.}\\
\item Another important question is whether there is an alternative approach to proving Conjecture~2 that avoids the notion of sparse decomposition, and even the notion of regular pairs. Such a programme has been  developed in the dense setting by Szemer\'edi and his collaborators, see~\cite{LeSaSz:Posa} for a particular instance of ``deregularizing'' a result originally resolved~\cite{KoSaSz:Posa} using the regularity method. However, this programme has not given a general alternative view, as of yet.
\end{itemize}
\section*{Acknowledgements}
The work on this project lasted from the beginning of 2008 and we are very grateful to the funding bodies for
their support. 

JH was funded by a BAYHOST fellowship, a DAAD fellowship, 
 Charles University grant GAUK 202-10/258009, EPSRC award EP/D063191/1, and by an EPSRC Postdoctoral Fellowship hosted by the Mathematics Institute at the University of Warwick. 
JK and ESz acknowledge the support of NSF grant
DMS-0902241.
DP acknowledges the support of the Marie Curie fellowship FIST,
 DFG grant TA 309/2-1, a DAAD fellowship,
 Czech Ministry of
Education project 1M0545,  EPSRC award EP/D063191/1, grant  PIEF-GA-2009-253925 of the European Union's FP7/2007-2013, and EPSRC
Additional Sponsorship EP/J501414/1.
MS was supported by a FAPESP fellowship, and by FAPESP travel grant  PQ-EX 2008/50338-0, also CMM-Basal, and  FONDECYT  grants 11090141 and 1140766.
\bibliographystyle{alpha}
\bibliography{LKSAnnou}

\newcommand{\etalchar}[1]{$^{#1}$}
\begin{thebibliography}{HKP{\etalchar{+}}d}

\bibitem[AKS95]{AKS95}
M.~Ajtai, J.~Koml{\'o}s, and E.~Szemer{\'e}di.
\newblock On a conjecture of {L}oebl.
\newblock In {\em Graph theory, combinatorics, and algorithms, Vol.\ 1, 2
  (Kalamazoo, MI, 1992)}, Wiley-Intersci. Publ., pages 1135--1146. Wiley, New
  York, 1995.

\bibitem[AKSS]{AKSS07+}
M.~Ajtai, J.~Koml\'os, M.~Simonovits, and E.~Szemer\'edi.
\newblock Proof of the {E}rd{\H{o}}s-{T}. {S}\'os conjecture for large trees.
\newblock In preparation.

\bibitem[ASS09]{AlShSt:DistinctyRegularityPartitions}
N.~Alon, A.~Shapira, and U.~Stav.
\newblock Can a graph have distinct regular partitions?
\newblock {\em SIAM J. Discrete Math.}, 23(1):278--287, 2008/09.

\bibitem[BCL09]{borgs-2008}
C.~Borgs, J.~Chayes, and L.~Lov\'asz.
\newblock Moments of two-variable functions and the uniqueness of graph limits.
\newblock {\em J. Geom. and Func. Anal}, 19:1597--1619, 2009.

\bibitem[BD96]{bradob}
S.~Brandt and E.~Dobson.
\newblock The {E}rd{\H o}s--{S}\'os conjecture for graphs of girth $5$.
\newblock {\em Discr. Math.}, 150:411--414, 1996.

\bibitem[BLW00]{BLW00}
C.~Bazgan, H.~Li, and M.~Wo{\'z}niak.
\newblock On the {L}oebl-{K}oml\'os-{S}\'os conjecture.
\newblock {\em J. Graph Theory}, 34(4):269--276, 2000.

\bibitem[Coo09]{Cooley08}
O.~Cooley.
\newblock Proof of the {L}oebl-{K}oml\'os-{S}\'os conjecture for large, dense
  graphs.
\newblock {\em Discrete Math.}, 309(21):6190--6228, 2009.

\bibitem[Dob02]{Dob02}
E.~Dobson.
\newblock Constructing trees in graphs whose complement has no {$K\sb {2,s}$}.
\newblock {\em Combin. Probab. Comput.}, 11(4):343--347, 2002.

\bibitem[EFLS95]{EFLS95}
P.~Erd{\H{o}}s, Z.~F{\"u}redi, M.~Loebl, and V.~T. S{\'o}s.
\newblock Discrepancy of trees.
\newblock {\em Studia Sci. Math. Hungar.}, 30(1-2):47--57, 1995.

\bibitem[EG59]{ErGa}
P.~Erd{\H{o}}s and T.~Gallai.
\newblock On maximal paths and circuits of graphs.
\newblock {\em Acta Math. Acad. Sci. Hungar}, 10:337--356 (unbound insert),
  1959.

\bibitem[FS13]{FurSim}
Z.~F{\"u}redi and M.~Simonovits.
\newblock The history of degenerate (bipartite) extremal graph problems.
\newblock In {\em Erd\"os centennial}, volume~25 of {\em Bolyai Soc. Math.
  Stud.}, pages 169--264. J\'anos Bolyai Math. Soc., Budapest, 2013.

\bibitem[GG67]{GerencserGyarfas}
L.~Gerencs{\'e}r and A.~Gy{\'a}rf{\'a}s.
\newblock On {R}amsey-type problems.
\newblock {\em Ann. Univ. Sci. Budapest. E\"otv\"os Sect. Math.}, 10:167--170,
  1967.

\bibitem[Hax01]{Haxell:TreeEmbeddings}
P.~E. Haxell.
\newblock Tree embeddings.
\newblock {\em J. Graph Theory}, 36(3):121--130, 2001.

\bibitem[HKP{\etalchar{+}}a]{LKSsparseapprox}
J.~Hladk{\'y}, J.~Koml{\'o}s, D.~Piguet, M.~Simonovits, M.~Stein, and
  E.~Szemer{\'e}di.
\newblock The approximate {L}oebl-{K}oml\'os-{S}\'os conjecture.
\newblock arXiv: 1211.3050.

\bibitem[HKP{\etalchar{+}}b]{LKSsparseapprox1}
J.~Hladk{\'y}, J.~Koml{\'o}s, D.~Piguet, M.~Simonovits, M.~Stein, and
  E.~Szemer{\'e}di.
\newblock The approximate {L}oebl-{K}oml\'os-{S}\'os conjecture {I}: {T}he
  sparse decomposition.
\newblock arXiv: 1408:3858.

\bibitem[HKP{\etalchar{+}}c]{LKSsparseapprox2}
J.~Hladk{\'y}, J.~Koml{\'o}s, D.~Piguet, M.~Simonovits, M.~Stein, and
  E.~Szemer{\'e}di.
\newblock The approximate {L}oebl-{K}oml\'os-{S}\'os conjecture {II}: {T}he
  rough structure of {LKS} graphs.
\newblock arXiv: 1408:3871.

\bibitem[HKP{\etalchar{+}}d]{LKSsparseapprox3}
J.~Hladk{\'y}, J.~Koml{\'o}s, D.~Piguet, M.~Simonovits, M.~Stein, and
  E.~Szemer{\'e}di.
\newblock The approximate {L}oebl-{K}oml\'os-{S}\'os conjecture {III}: {T}he
  finer structure of {LKS} graphs.
\newblock arXiv: 1408:3866.

\bibitem[HKP{\etalchar{+}}e]{LKSsparseapprox4}
J.~Hladk{\'y}, J.~Koml{\'o}s, D.~Piguet, M.~Simonovits, M.~Stein, and
  E.~Szemer{\'e}di.
\newblock The approximate {L}oebl-{K}oml\'os-{S}\'os conjecture {IV}:
  {E}mbedding techniques and the proof of the main result.
\newblock arXiv: 1408:3870.

\bibitem[HLT02]{HLT02}
P.~E. Haxell, T.~Luczak, and P.~W. Tingley.
\newblock Ramsey numbers for trees of small maximum degree.
\newblock {\em Combinatorica}, 22(2):287--320, 2002.
\newblock Special issue: Paul Erd\H os and his mathematics.

\bibitem[HP]{HladkyPiguet}
J.~Hladk{\'y} and D.~Piguet.
\newblock {L}oebl-{K}oml\'os-{S}\'os {C}onjecture: dense case.
\newblock arXiv:0805.4834.

\bibitem[KO09]{KuOs}
D.~K{\"u}hn and D.~Osthus.
\newblock Embedding large subgraphs into dense graphs.
\newblock In {\em Surveys in combinatorics 2009}, volume 365 of {\em London
  Math. Soc. Lecture Note Ser.}, pages 137--167. Cambridge Univ. Press,
  Cambridge, 2009.

\bibitem[KSS98a]{PoSe}
J.~Koml{\'o}s, G.~S{\'a}rk{\"o}zy, and E.~Szemer{\'e}di.
\newblock Proof of the {S}eymour conjecture for large graphs.
\newblock {\em Ann. Comb.}, 2(1):43--60, 1998.

\bibitem[KSS98b]{KoSaSz:Posa}
J.~Koml{\'o}s, G.~N. S{\'a}rk{\"o}zy, and E.~Szemer{\'e}di.
\newblock Proof of the {S}eymour conjecture for large graphs.
\newblock {\em Ann. Comb.}, 2(1):43--60, 1998.

\bibitem[KSSS02]{regu2}
J.~Koml{\'o}s, A.~Shokoufandeh, M.~Simonovits, and E.~Szemer{\'e}di.
\newblock The regularity lemma and its applications in graph theory.
\newblock In {\em Theoretical aspects of computer science (Tehran, 2000)},
  volume 2292 of {\em Lecture Notes in Comput. Sci.}, pages 84--112. Springer,
  Berlin, 2002.

\bibitem[LSS10]{LeSaSz:Posa}
I.~Levitt, G.~N. S{\'a}rk{\"o}zy, and E.~Szemer{\'e}di.
\newblock How to avoid using the regularity lemma: {P}\'osa's conjecture
  revisited.
\newblock {\em Discrete Math.}, 310(3):630--641, 2010.

\bibitem[PS08]{PS2}
D.~Piguet and M.~J. Stein.
\newblock {L}oebl-{K}oml\'os-{S}\'os conjecture for trees of diameter 5.
\newblock {\em Electron. J. Combin.}, 15(1):Research Paper 106, 11 pp.
  (electronic), 2008.

\bibitem[PS12]{PS07+}
D.~Piguet and M.~J. Stein.
\newblock An approximate version of the {L}oebl-{K}oml\'os-{S}\'os conjecture.
\newblock {\em J. Combin. Theory Ser. B}, 102(1):102--125, 2012.

\bibitem[Sof00]{Sof00}
S.~N. Soffer.
\newblock The {K}oml\'os-{S}\'os conjecture for graphs of girth 7.
\newblock {\em Discrete Math.}, 214(1--3):279--283, 2000.

\bibitem[SW97]{sacwoz}
J.-F. Sacl\'e and M.~Wo\'zniak.
\newblock A note on the {E}rd{\H o}s--{S}\'os conjecture for graphs without
  ${C}_4$.
\newblock {\em J.~Combin.\ Theory (Series B)}, 70(2):229--234, 1997.

\bibitem[Sze78]{Sze78}
E.~Szemer{\'e}di.
\newblock Regular partitions of graphs.
\newblock In {\em Probl\`emes combinatoires et th\'eorie des graphes (Colloq.
  Internat. CNRS, Univ. Orsay, Orsay, 1976)}, volume 260 of {\em Colloq.
  Internat. CNRS}, pages 399--401. CNRS, Paris, 1978.

\bibitem[Wo{\'z}96]{woz}
M.~Wo{\'z}niak.
\newblock On the {E}rd{\H o}s--{S}\'os conjecture.
\newblock {\em J.~Graph Theory}, 21(2):229--234, 1996.

\bibitem[Zha11]{Z07+}
Y.~Zhao.
\newblock Proof of the {$(n/2-n/2-n/2)$} conjecture for large {$n$}.
\newblock {\em Electron. J. Combin.}, 18(1):Paper 27, 61, 2011.

\end{thebibliography}

\end{document}